\documentclass[english]{article}
\usepackage[T1]{fontenc}
\usepackage[latin9]{inputenc}
\usepackage{geometry}
\usepackage{float}
\usepackage{amsmath}
\usepackage{amsthm}
\usepackage{amssymb}
\usepackage{graphicx}

\makeatletter
\numberwithin{equation}{section}
\numberwithin{figure}{section}
\theoremstyle{plain}
\newtheorem{thm}{\protect\theoremname}
\theoremstyle{plain}
\newtheorem{prop}[thm]{\protect\propositionname}
\theoremstyle{remark}
\newtheorem{rem}[thm]{\protect\remarkname}

\allowdisplaybreaks
\usepackage{mathabx}

\makeatother

\usepackage{babel}
\providecommand{\propositionname}{Proposition}
\providecommand{\remarkname}{Remark}
\providecommand{\theoremname}{Theorem}

\begin{document}
\global\long\def\ve{\varepsilon}%
\global\long\def\R{\mathbb{R}}%
\global\long\def\C{\mathbb{C}}%
\global\long\def\F{\mathbb{F}}%
\global\long\def\Rn{\mathbb{R}^{n}}%
\global\long\def\Rd{\mathbb{R}^{d}}%
\global\long\def\E{\mathbb{E}}%
\global\long\def\P{\mathbb{P}}%
\global\long\def\bx{\mathbf{x}}%
\global\long\def\vp{\varphi}%
\global\long\def\ra{\rightarrow}%
\global\long\def\smooth{C^{\infty}}%
\global\long\def\Tr{\mathrm{Tr}}%
\global\long\def\bra#1{\left\langle #1\right|}%
\global\long\def\ket#1{\left|#1\right\rangle }%
\global\long\def\ud{\mathrm{d}}%
\global\long\def\blam{\boldsymbol{\lambda}}%

\title{Fast randomized entropically regularized\\
 semidefinite programming}
\author{Michael Lindsey\\
{\footnotesize{}UC Berkeley}\\
{\footnotesize{}$\texttt{lindsey@berkeley.edu}$}}
\maketitle
\begin{abstract}
We develop a practical approach to semidefinite programming (SDP)
that includes the von Neumann entropy, or an appropriate variant,
as a regularization term. In particular we solve the dual of the regularized
program, demonstrating how a carefully chosen randomized trace estimator
can be used to estimate dual gradients effectively. We also introduce
specialized optimization approaches for common SDP, specifically SDP
with diagonal constraint and the problem of the determining the spectral
projector onto the span of extremal eigenvectors. We validate our
approach on such problems with applications to combinatorial optimization
and spectral embedding.
\end{abstract}

\section{Introduction\label{sec:Introduction}}

In this work we are interested in the solution of semidefinite programs
(SDPs) of the form 
\begin{align*}
\underset{X\in\mathbb{R}^{n\times n}}{\text{minimize}}\quad & \Tr[CX]\\
\text{subject to}\quad & X\succeq0,\ \Tr\left[A_{k}X\right]=b_{k},\ k=1,\ldots,m.
\end{align*}
 A major obstacle in semidefinite programming is that enforcing the
semidefinite constraint on an $n\times n$ matrix $X$ in general
requires $\sim n^{3}$ operations. This can be understood intuitively
by realizing that detecting semidefiniteness (i.e., checking whether
all eigenvalues are nonnegative) essentially involves computing a
complete spectrum.

Major progress in avoiding $n^{3}$ scaling has been made in the case
when the optimal solution $X^{\star}$ is of rank $k\ll n$. In this
case, many specialized approaches are able to achieve per-iteration
complexity of $O(nk^{2})$. These include the Burer-Monteiro / SDPLR
method \cite{BurerMonteiro} and alternative manifold-constrained
optimization approaches \cite{manopt}, as well as the randomized
approach SketchyCGAL \cite{SketchyCGAL}.

In this work we are mainly interested in settings where the numerical
rank is high enough such that $O(nk^{2})$ scaling is prohibitive,
and we pursue algorithms that do not depend on any rank truncation.
However, our regularization approach also smooths the problem and
could possibly yield a convergence advantage even for low-rank problems,
provided that approximation is acceptable.

Inspired by the success of entropic regularization \cite{Cuturi2013}
of the Kantorovich problem of optimal transport, which is a linear
program, we pursue the regularization of SDPs by the von Neumann entropy,
as well as other variants as appropriate. Regularization by the von
Neumann entropy has been considered in \cite{LinLindsey2020,GibbsManifolds,Pavlov2023}
and moreover can be viewed as fundamental to the perspective of quantum
statistical mechanics at finite temperature (cf. \cite{Dissertation}
for mathematical introduction). We will develop an appropriate theory
of duality which shall be of no surprise from a physicist's point
of view, but which to our knowledge has not been exploited for the
purpose of fast solution of SDPs.

In contrast to the situation of entropically regularized optimal transport,
the duality theory of entropically regularized SDPs does not immediately
suggest any algorithmic treatment. Indeed, the dual gradients require
the computation of certain matrix traces that cannot be computed exactly
in less than $O(n^{3})$ time.

The idea for trace estimation is based on Hutchinson's trace estimator
\cite{Hutchinson,Hutch++}, with a `square root' trick that has also
been applied in the context of Gaussian process regression (GPR) \cite{SqrtHutch}.
Unlike the context of GPR in which the trick requires expensive matrix-vector
multiplications by the square root of the kernel matrix, in our context
the square root trick imposes no additional computational burden relative
to the `default' strategy.

We focus on two problem types of interest. The first type is that
of SDPs with diagonal constraint, which include the Goemans-Williamson
relaxation of the Max-Cut problem \cite{GoemansWilliamson}. For SDPs
with this structure we introduce a specialized optimization approach
that takes some loose inspiration from matrix scaling \cite{Cuturi2013}
and can be viewed as solving a sequence of minorized dual problems.

The second type is the SDP formulation of the problem of computing
the spectral projector onto the span of extremal eigenvectors. Here
we use the binary von Neumann entropy as a regularizer, which connects
to the theory of single-particle fermionic statistical mechanics (again
cf. \cite{Dissertation} for mathematical introduction), and our specialized
optimization approach is simply Newton's method, with updates computed
by appropriate randomized trace estimation. We demonstrate an application
to graph spectral embedding \cite{PothenEtAl1990,SpectralClustering}.

For fixed value of the regularization parameter, the per-iteration
cost of our algorithms scales according to the cost of matrix-vector
multiplication by the cost and constraint matrices $C$ and $A_{1},\ldots,A_{m}$,
yielding $O(n)$ per-iteration scaling in our applications for graphs
of bounded degree. Empirically we find that for the problem families
of interest, when the regularization parameter is fixed, the optimization
converges in $O(1)$ iterations.

We conclude the introduction with an outline of the paper. In Section
\ref{sec:entropicreg} we introduce the general framework of entropic
regularization and duality. In Section \ref{sec:trace} we describe
the trace estimation procedure needed in our algorithms, which is
equipped with rigorous concentration bounds. In Section \ref{sec:specialized_opt}
we describe specialized optimization approaches for problems of interest.
In Section \ref{sec:Applications} we discuss applications to the
Max-Cut problem and spectral clustering and present numerical experiments.

\section{Entropic regularization of SDP \label{sec:entropicreg}}

Consider the general semidefinite program 
\begin{align}
\underset{X\in\mathbb{R}^{n\times n}}{\text{minimize}}\quad & \Tr[CX]\label{eq:sdp}\\
\text{subject to}\quad & X\succeq0,\nonumber \\
 & \Tr\left[\mathbf{A}X\right]=\mathbf{b}.\nonumber 
\end{align}
Here $\mathbf{A}=(A_{1},\ldots,A_{m})\in\R^{m\times(n\times n)}$
indicates a vector of matrices, which together with $\mathbf{b}=(b_{1},\ldots,b_{m})\in\R^{m}$,
specifies $m$ linear equality constraints $\Tr[A_{i}X]=b_{i}$, $i=1,\ldots,m$,
on the optimization variable $X$. We can always assume that the cost
matrix $C\in\R^{n\times n}$ is symmetric by symmetrizing if necessary,
which does not alter the objective. We moreover assume without loss
of generality that $A_{1},\ldots,A_{m}$ are symmetric as they too
can be symmetrized without altering the expression $\Tr[\mathbf{A}X]$.

Algorithmically we shall only require access to the cost matrix $C$
and the constraint matrices $A_{i}$ via matrix-vector multiplications
(matvecs). Thus in the case where the cost and constraint matrices
are sparse or otherwise structured, fast matvecs can be easily exploited.

More generally we can consider complex-valued $C,\mathbf{A},\mathbf{b}$
and Hermitian positive definite optimization variable $X$, where
$C$ can be assumed to be Hermitian, and the objective is replaced
by $\mathrm{Re}\left(\Tr[CX]\right)$. For simplicity we restrict
our discussion to the real case.

Consider regularizing the problem by the addition of a von Neumann
entropy term 
\[
S(X):=\Tr\left[X\log X-X\right]
\]
 as follows: 
\begin{align*}
\underset{X\in\mathbb{R}^{n\times n}}{\text{minimize}}\quad & \Tr[CX]+\beta^{-1}S(X)\\
\text{subject to}\quad & \Tr\left[\mathbf{A}X\right]=\mathbf{b}.
\end{align*}
 Here $\beta\in(0,\infty)$ carries the physical interpretation of
an inverse temperature in quantum statistical mechanics, and the domain
of $S(X)$ is implicitly understood to be the positive definite cone
$\{X\succ0\}$.

\subsection{Duality}

Introducing a Lagrange multiplier $\boldsymbol{\lambda}\in\R^{m}$,
we obtain the minimax problem for the suitable Lagrangian 
\[
\min_{X\in\R^{n\times n}}\max_{\blam\in\R^{m}}\mathcal{L}(X,\blam),
\]
 where 
\begin{align*}
\mathcal{L}(X,\blam):=\  & \Tr[CX]+\beta^{-1}S(X)+\blam\cdot\left(\mathbf{b}-\Tr[\mathbf{A}\cdot X]\right)\\
=\  & \Tr\left[\left(C-\blam\cdot\mathbf{A}\right)X\right]+\beta^{-1}S(X)+\mathbf{b}\cdot\blam.
\end{align*}

Minimization over $X$ yields optimizer 
\[
X=e^{-\beta(C-\blam\cdot\mathbf{A})},
\]
 and the dual objective $g(\blam)$ is then computed as 
\[
g(\blam):=\max_{X}\mathcal{L}(X,\blam)=\mathbf{b}\cdot\blam-\beta^{-1}\Tr\left[e^{-\beta(C-\blam\cdot\mathbf{A})}\right],
\]
 yielding the dual problem 
\begin{align*}
\underset{\blam\in\R^{m}}{\text{maximize}}\quad & \mathbf{b}\cdot\blam-\beta^{-1}\Tr\left[e^{-\beta(C-\blam\cdot\mathbf{A})}\right],
\end{align*}
 in which the second term of the objective carries the interpretation
of the quantum Gibbs free energy.

When a dual solution $\blam^{\star}$ is obtained, the primal solution
$X^{\star}=X(\blam^{\star})$, where 
\begin{equation}
X(\blam):=e^{-\beta(C-\blam\cdot\mathbf{A})}\label{eq:Xblam}
\end{equation}
 carries the interpretation of an unnormalized quantum density operator.

\subsection{Dual gradients}

To solve the dual optimization problem via first-order methods, we
are motivated to compute the gradient of $g(\blam)$. These can be
obtained analytically as 
\begin{equation}
\nabla g(\blam)=\nabla_{\blam}\mathcal{L}(X(\blam),\blam)=\mathbf{b}-\Tr\left[\mathbf{A}X(\blam)\right].\label{eq:grad}
\end{equation}
 It is inefficient to evaluate the gradient directly, since forming
the matrix exponential exactly requires $O(n^{3})$ operations. We
will discuss in Section \ref{sec:trace} how to estimate traces of
the form appearing in (\ref{eq:grad}) using only matvecs by the cost
and constraint matrices.

\subsection{The case of diagonal constraint}

Of particular interest is the case where $A_{i}=e_{i}e_{i}^{\top}$,
$i=1,\ldots,n$, i.e., $m=n$. In the case the primal SDP takes the
form 
\begin{align}
\underset{X\in\mathbb{R}^{n\times n}}{\text{minimize}}\quad & \Tr[CX]\label{eq:sdpdiag}\\
\text{subject to}\quad & X\succeq0,\nonumber \\
 & \mathrm{diag}(X)=\mathbf{b},\nonumber 
\end{align}
 the dual regularized SDP takes the form 
\[
\underset{\blam\in\R^{m}}{\text{maximize}}\quad\mathbf{b}\cdot\blam-\beta^{-1}\Tr\left[e^{-\beta(C-\mathrm{diag}(\blam))}\right],
\]
 and the dual gradients are 
\[
\nabla g(\blam)=\mathbf{b}-\mathrm{diag}\left[X(\blam)\right].
\]
 Hence computing the dual gradients requires the estimation of the
diagonal of a positive definite matrix, which is in particular presented
as a matrix exponential.

\subsection{Binary-entropic regularization}

Although (\ref{eq:sdp}) is the general form of an SDP, many specific
SDPs can only be reduced to this form via the introduction of additional
optimization variables and a mess of dense equality constraints. We
provide an example of how the framework of entropically regularized
SDP can extend to another setting of interest without sacrificing
conceptual clarity and algorithmic efficiency. Consider the SDP 

\begin{align}
\underset{X\in\mathbb{R}^{n\times n}}{\text{minimize}}\quad & \Tr[CX]\label{eq:sdpBin}\\
\text{subject to}\quad & 0\preceq X\preceq I,\nonumber \\
 & \Tr\left[\mathbf{A}X\right]=\mathbf{b},\nonumber 
\end{align}
 which differs only from (\ref{eq:sdp}) via the inclusion of an upper
bound on $X$ in the Loewner ordering.

In this setting, we may instead regularize by the \emph{binary }von
Neumann entropy, which generalizes the ordinary binary entropy to
the matrix case:
\[
S_{\mathrm{bin}}(X):=\Tr[X\log X]+\Tr\log[(I-X)\log(I-X)].
\]
 The domain of $S_{\mathrm{bin}}$ is understood to be $\{0\prec X\prec I\}$,
the set of all symmetric matrices with eigenvalues strictly between
$0$ and $1$.

In this setting the Lagrangian reads as 
\[
\mathcal{L}(X,\blam)=\ =\ \Tr\left[\left(C-\blam\cdot\mathbf{A}\right)X\right]+\beta^{-1}S_{\mathrm{bin}}(X)+\mathbf{b}\cdot\blam,
\]
 whose primal optimizer for fixed $\blam$ is 
\begin{equation}
X(\blam)=\left[I+e^{\beta\left(C-\blam\cdot\mathbf{A}\right)}\right]^{-1}=F_{\beta}\left(C-\blam\cdot\mathbf{A}\right),\label{eq:XblamBin}
\end{equation}
 where 
\[
F_{\beta}(x)=\frac{1}{1+e^{\beta x}}=\frac{1}{2}\left(1-\tanh(\beta x/2)\right)
\]
 is the Fermi-Dirac function at inverse temperature $\beta$.

Some computation reveals that the dual objective can be written 
\[
g(\blam)=\mathbf{b}\cdot\blam+\beta^{-1}\Tr\left[\log\left(I+e^{-\beta(C-\blam\cdot\mathbf{A})}\right)\right],
\]
 in which the last term can be interpreted as a fermionic free energy.
The dual gradient is once again given by the expression 

\[
\nabla g(\blam)=\nabla_{\blam}\mathcal{L}(X(\blam),\blam)=\mathbf{b}-\Tr\left[\mathbf{A}X(\blam)\right],
\]
 where now $X(\blam)$ carries a different interpretation, after (\ref{eq:XblamBin}).

.

\subsection{Extremal eigenvalue problem}

A special case of (\ref{eq:sdpBin}) of particular interest is the
\begin{align}
\underset{X\in\mathbb{R}^{n\times n}}{\text{minimize}}\quad & \Tr[CX]\label{eq:sdpeig}\\
\text{subject to}\quad & 0\preceq X\preceq I,\nonumber \\
 & \Tr\left[X\right]=k,\nonumber 
\end{align}
 which is an SDP formulation of the problem of finding the lowest
$k$ eigenvectors of a symmetric matrix $C$. Indeed, the optimal
solution is recovered as 
\[
X^{\star}=\sum_{i=1}^{k}u_{i}u_{i}^{\top},
\]
 where $(\lambda_{i},u_{i})$, $i=1,\ldots,n$, are the eigenpairs
of $C$, ordered $\lambda_{1}\leq\cdots\leq\lambda_{n}$, assuming
we have a gap $\lambda_{k}<\lambda_{k+1}$ between the $k$-th and
$(k+1)$-th eigenvalues.

The format (\ref{eq:sdpBin}) is recovered by taking $m=1$, $A_{1}=I$,
and $b_{1}=k$. In this case the Lagrange multiplier is a scalar carrying
the interpretation of a chemical potential, hence we shall denote
it by $\mu$, and the dual objective can be written 
\[
g(\mu)=k\mu+\beta^{-1}\Tr\left[\log\left(1+e^{-\beta(C-\mu I)}\right)\right],
\]
 and the derivative is 
\begin{equation}
g'(\mu)=k-\Tr\left[X(\mu)\right].\label{eq:eiggrad}
\end{equation}
 In fact in this case, the second derivative admits a simple expression,
owing to the fact that the constraint matrix $I$ commutes with the
cost matrix $C$: 
\[
g''(\mu)=-\beta\,\Tr\left[\frac{e^{\beta(C-\mu I)}}{\left[I+e^{\beta(C-\mu I)}\right]^{2}}\right],
\]
 which can be rewritten conveniently as 
\begin{equation}
g''(\mu)=-\beta\,\Tr\left[\left(I-X(\mu)\right)X(\mu)\right],\label{eq:eighess}
\end{equation}
 where we recall $X(\mu)=\left[I+e^{\beta(C-\mu I)}\right]^{-1}$.

\section{Trace estimation \label{sec:trace}}

It is inefficient to exactly construct the entire matrix $X(\blam)$
following (\ref{eq:sdp}). We could instead pursue a randomized approach
for computing the gradient (\ref{eq:grad}) directly without forming
$X(\blam)$.

Instead, we will only ever require the capacity to perform matvecs
by $Y(\blam):=X^{1/2}(\blam)$. In fact, we only require the capacity
to perform matvecs by any $Y(\blam)$ achieving the factorization
$X(\blam)=Y(\blam)Y(\blam)^{\top}$. However, given the presentation
of $X(\blam)$ as a matrix function, the most straightforward practical
choice of factorization available is furnished by $Y(\blam)=X^{1/2}(\blam)$,
and we will restrict our discussion to this choice.

In the case (\ref{eq:Xblam}) of ordinary entropic regularization
we have 
\begin{equation}
Y(\blam)=e^{-\frac{\beta}{2}(C-\blam\cdot\mathbf{A})},\label{eq:Y}
\end{equation}
 while in the case (\ref{eq:XblamBin}) of binary entropic regularization
we have 
\begin{equation}
Y(\blam)=F_{\beta}^{1/2}(C-\blam\cdot\mathbf{A}),\label{eq:Ybin}
\end{equation}
 where $F_{\beta}^{1/2}$ denotes the square root of the Fermi-Dirac
function 
\begin{equation}
F_{\beta}^{1/2}(x)=\frac{1}{\sqrt{1+e^{\beta x}}}.\label{eq:sqrtfermidirac}
\end{equation}
Note that in either case, $Y(\blam)$ presents as a matrix function
of $C-\blam\cdot\mathbf{A}$, which is no more difficult to deal with
computationally than $X(\blam)$ itself.

The problem of estimating the diagonal of a matrix $X$ using only
matvecs by $X$ has been considered in \cite{BekasSaad2007}, for
example, using a Hutchinson-type estimator. The variance of this estimator,
however, depends on the locality of the matrix whose diagonal is to
be estimated. By using a Hutchinson-type estimator that instead relies
on matvecs by the \emph{square root} matrix, the variance is guaranteed
to be improved \cite{SqrtHutch} for all traces $\Tr[\mathbf{A}X]$.
In the case of diagonal estimation, as we shall see, the diagonal
entries can in fact be recovered with relative variance that is \emph{universal}.
In the setting of \cite{SqrtHutch}, one downside of the square root
approach is that matvecs by the square root are more expensive than
matvecs by $X$. In our setting, there is essentially no difference
in cost as both must be constructed using similar matrix function
techniques.

\subsection{Estimator}

To derive our estimator, compute 
\begin{align*}
\Tr[\mathbf{A}X(\blam)] & =\Tr[Y(\blam)\mathbf{A}Y(\blam)]\\
 & =\E_{z\sim\mathcal{N}(0,I)}\left(z^{\top}Y(\blam)\mathbf{A}Y(\blam)z\right)\\
 & =\E_{z\sim\mathcal{N}(0,I)}\left(\left[Y(\blam)z\right]^{\top}\mathbf{A}\left[Y(\blam)z\right]\right).
\end{align*}
This equation suggests the unbiased estimator 
\begin{equation}
\mathbf{a}^{(N)}(\blam):=\frac{1}{N}\sum_{i=1}^{N}\left[Y(\blam)z^{(i)}\right]^{\top}\mathbf{A}\left[Y(\blam)z^{(i)}\right]\label{eq:est}
\end{equation}
 for $\Tr[\mathbf{A}X(\blam)]$, where $z^{(i)}$, $i=1,\ldots,N$,
are independent and identically distributed according to the standard
normal distribution $\mathcal{N}(0,I)$.

It is possible more generally to adopt as the distribution for the
$z^{(i)}$ any distribution yielding independent entries of unit variance,
such as the Rademacher distribution taking values in $\{\pm1\}$.
For simplicity, we restrict our practical attention to the Gaussian
case.

Note for concreteness that in the case of the diagonal constraint
(where $m=n$ and $A_{k}=e_{k}e_{k}^{\top}$, $k=1,\ldots,n$), we
can write 
\[
\mathbf{a}^{(N)}(\blam)=\frac{1}{N}\sum_{i=1}^{N}\left[Y(\blam)z^{(i)}\right]^{\odot2},
\]
 where $u^{\odot2}$ indicates the entrywise second power of $u$.
If we let $Z=[z^{(1)}\,\cdots\,z^{(N)}]\in\R^{n\times N}$, we can
more compactly write 
\begin{equation}
\mathbf{a}^{(N)}(\blam)=\frac{1}{N}\left[Y(\blam)Z\right]^{\odot2}\,\mathbf{1}_{N},\label{eq:diagest}
\end{equation}
 where $\mathbf{1}_{N}\in\R^{N}$ is the vector of all ones.

\subsection{Covariance and concentration inequalities}

In this section let us fix $X=X(\blam)$, $Y=Y(\blam)$, and $\mathbf{a}^{(N)}=\mathbf{a}^{(N)}(\blam)$
for notational simplicity.

For $z\sim N(0,I)$, we have 
\[
\E[(Yz)^{\top}\mathbf{A}(Yz)]=\Tr[\mathbf{A}X],
\]
 and moreover it can be readily verified via Wick's theorem, which
yields the identity $\E[z_{i}z_{j}z_{k}z_{l}]=\delta_{ij}\delta_{kl}+\delta_{ik}\delta_{il}+\delta_{il}\delta_{jk}$,
that the covariance of $(Yz)^{\top}\mathbf{A}(Yz)\in\R^{m}$ is given
by 
\[
\Sigma:=\mathrm{Cov}[(Yz)^{\top}\mathbf{A}(Yz)]=\left(2\,\Tr[A_{i}XA_{j}X]\right)_{i,j=1,\ldots,m}.
\]
Therefore the mean and covariance of our estimator $\mathbf{a}^{(N)}$
can be recovered as 
\[
\E[\mathbf{a}^{(N)}]=\Tr[\mathbf{A}X],\quad\mathrm{Cov}(\mathbf{a}^{(N)})=N^{-1}\Sigma.
\]

We moreover have the following concentration bounds on the estimator,
applicable for more general choice of the distribution of the $z^{(i)}$: 
\begin{prop}
Let $\delta\in(0,1/2]$, and suppose that $z^{(i)},i=1,\ldots,N$
are i.i.d. with entries that are themselves i.i.d. sub-Gaussian random
variables with a constant sub-Gaussian parameter, mean zero, and unit
variance. There exist constants $c,C>0$ such that if $N>c\log\left(\frac{m}{\delta}\right)$,
then with probability at least $1-\delta$, the inequality 
\[
\vert a_{i}^{(N)}(\blam)-\Tr\left[A_{i}X(\blam)\right]\vert\leq C\sqrt{\log\left(\frac{m}{\delta}\right)}\ \sqrt{\frac{\Tr[A_{i}XA_{i}X]}{N}}
\]
 holds for all $i=1,\ldots,m$.

In the case $\mathbf{A}=(e_{1}e_{1}^{\top},\ldots,e_{n}e_{n}^{\top})$
yielding the SDP with diagonal constraint (\ref{eq:sdpdiag}), it
follows that under the same conditions, 
\[
\frac{\left\Vert \mathbf{a}^{(N)}(\blam)-\mathrm{diag}\left[X(\blam)\right]\right\Vert _{2}}{\left\Vert \mathrm{diag}\left[X(\blam)\right]\right\Vert _{2}}\leq C\sqrt{\frac{\log\left(\frac{n}{\delta}\right)}{N}}.
\]
\end{prop}

\begin{rem}
Note that the second statement provides a bound on the \emph{relative}
error of the diagonal estimator that is \emph{universal} apart from
logarithmic dependence on $n$. Unlike the estimator of \cite{BekasSaad2007},
the error does not depend on locality of $X$.
\end{rem}

\begin{proof}
The first statement results from direct application of Lemma 2.1 of
\cite{Hutch++} to the Hutchinson estimator for each of $YA_{i}Y$,
$i=1,\ldots,m$, together with the union bound over the $m$ separate
traces.

Note that the first statement implies that 
\[
\left\Vert \mathbf{a}^{(N)}(\blam)-\mathrm{diag}\left[X(\blam)\right]\right\Vert _{2}\leq C\sqrt{\log\left(\frac{m}{\delta}\right)}\ \sqrt{\frac{\sum_{i=1}^{m}\Tr[A_{i}XA_{i}X]}{N}},
\]
 and the second statement follows immediately in the case $\mathbf{A}=(e_{1}e_{1}^{\top},\ldots,e_{n}e_{n}^{\top})$.
\end{proof}

\subsection{Matvecs \label{sec:matvecs}}

Several pieces remain in the specification of a concrete algorithm
for entropically regularized semidefinite programming. First we must
explain how to perform matvecs by $Y(\blam)$.

Observe that in the basic case (\ref{eq:Y}), $Y(\blam)=e^{-\beta(C-\blam\cdot\mathbf{A})}$
is defined as a matrix exponential of a matrix $C-\blam\cdot\mathbf{A}$
by which matvecs are assumed to be efficient, since they reduce to
matvecs by the cost and constraint matrices. In many cases of interest,
$C-\blam\cdot\mathbf{A}$ is in fact a sparse matrix with $O(n)$
nonzero entries. Therefore we apply the algorithm of \cite{expmv}
which lifts a matvec routine for a matrix to a matvec routine for
its exponential. This algorithm is based on splitting the exponential
into a product of exponentials closer to the identity, which are in
turn approximated by a Taylor series.

In the case of (\ref{eq:Ybin}), $Y(\blam)=F_{\beta}^{1/2}(C-\blam\cdot\mathbf{A})$
is furnished as a more complicated matrix function of the same matrix.
Recall that $F_{\beta}^{1/2}$ is the square root of the Fermi-Dirac
function, defined as in (\ref{eq:sqrtfermidirac}). For this matrix
function, no specialized algorithm such as \cite{expmv} is available,
so given lower and upper bounds on the spectrum of $C-\blam\cdot\mathbf{A}$,
we simply approximate $F_{\beta}^{1/2}$ via minimax polynomial approximation
\cite{Chebfun} to desired tolerance on the appropriate interval,
yielding $F_{\beta}^{1/2}\approx\sum_{k=0}^{K}a_{k}p_{k}$ on this
interval where $p_{k}$ are the appropriately shifted and scaled Chebyshev
polynomials. Then matvecs by $Y(\blam)\,u$ can be achieved by constructing
the matvecs $v_{k}=p_{k}(C-\blam\cdot\mathbf{A})\,u$ using the three-term
recurrence and linearly combining the results as $\sum_{k=0}^{K}a_{k}v_{k}$.

In the case of the extremal eigenvalue problem (\ref{eq:sdpeig}),
an interval $[B_{1},B_{2}]$ bounding the spectrum of the cost matrix
$C$, the spectrum of $C-\mu I$ is contained within $[B_{1}-\mu,B_{2}-\mu]$.
Moreover, the optimal $\mu^{\star}$ must lie in $[B_{1},B_{2}]$,
so $B_{1}-\mu^{\star}\geq B_{1}-B_{2}$, and $B_{2}-\mu^{\star}\leq B_{2}-B_{1}$.
Thus near the optimizer $\mu^{\star}$ we can use the interval $[-R,R]$
for polynomial approximation, where $R=B_{2}-B_{1}>0$ is a bound
on the spectral range.

\section{Specialized optimization algorithms \label{sec:specialized_opt}}

The trace estimator (\ref{eq:est}) introduced in Section \ref{sec:trace}
allows us to compute unbiased estimators of dual gradients following
(\ref{eq:grad}), enabling stochastic first-order optimization methods
such as stochastic gradient descent (SGD) and its accelerated variants.
An exploration of such approaches in a general setting is reserved
for future work. In this section we discuss how in two settings of
particular interest, specialized optimization approaches making use
of trace estimators are available.

\subsection{Diagonal constraint: noncommutative matrix scaling \label{sec:scaling}}

In the case of the SDP with diagonal constraint (\ref{eq:sdpdiag}),
recall that we seek $\blam$ such that 
\[
\mathrm{diag}\left[X(\blam)\right]=\mathrm{diag}\left[e^{-\beta(C-\mathrm{diag}(\blam))}\right]=\mathbf{b}.
\]
 One potential approach takes some inspiration from the idea of matrix
scaling algorithms such as Sinkhorn scaling \cite{Cuturi2013}. Given
some guess $\blam^{0}$, observe that there exists a unique perturbation
$\Delta\blam$ such that 
\[
\mathrm{diag}\left[e^{\beta\mathrm{diag}(\Delta\blam)}e^{-\beta(C-\mathrm{diag}(\blam^{0}))}\right]=\mathbf{b},
\]
 which is given by 
\begin{equation}
\Delta\blam=\beta^{-1}\left(\log\mathbf{b}-\log\mathrm{diag}\left[X(\blam^{0})\right]\right),\label{eq:deltablam}
\end{equation}
 where $\log$ indicates the entrywise logarithm on vector input.

If $C$ were diagonal (hence commuting with all diagonal matrices),
we could replace
\[
e^{\beta\mathrm{diag}(\Delta\blam)}e^{-\beta(C-\mathrm{diag}(\blam^{0}))}\ra e^{-\beta(C-\mathrm{diag}(\blam^{0}+\Delta\blam))},
\]
 i.e., 
\[
e^{\beta\mathrm{diag}(\Delta\blam)}X(\blam^{0})\ra X(\blam^{0}+\Delta\blam),
\]
 which suggests the update 
\begin{equation}
\blam^{1}\leftarrow\blam^{0}+\Delta\blam\label{eq:blam1}
\end{equation}
 that achieves the diagonal constraint $\mathrm{diag}[X(\blam^{1})]=\mathbf{b}$
exactly.

In reality, $C$ is not diagonal, but we can still use this update.
Interestingly, a more rigorous theoretical justification for this
update \emph{does not }follow from any argument based Suzuki-Trotter
expansion, since the matrix $C-\mathrm{diag}(\blam^{0})$ is not of
small size even when the guess $\blam^{0}$ is nearly optimal.

Instead, the update rule, which we call \emph{noncommutative matrix
scaling}, can be interpreted as exactly solving a \emph{minorization}
of the dual maximization problem about the guess $\blam^{0}$.

Indeed, recall the dual objective 
\[
g(\blam)=\mathbf{b}\cdot\blam-\beta^{-1}\Tr\left[e^{-\beta(C-\mathrm{diag}(\blam))}\right].
\]
 For a given guess $\blam^{0}$, we will define $g^{0}$ satisfying
the following properties: 
\begin{enumerate}
\item $g^{0}$ is strictly concave;
\item $g^{0}\leq g$ pointwise, with $g^{0}(\blam^{0})=g(\blam)$; 
\item With $\blam^{1}$ defined in terms of $\blam^{0}$ following (\ref{eq:blam1})
and (\ref{eq:deltablam}), we have 
\begin{equation}
\blam^{1}=\underset{\blam}{\text{argmax}}\ g^{0}(\blam);\label{eq:minor}
\end{equation}
\item $\blam^{1}=\blam^{0}$ if and only if $\blam^{0}=\blam^{\star}$.
\end{enumerate}
From these conditions on $g^{0}$, it follows that 
\[
g(\blam^{1})\geq g^{0}(\blam^{1})\geq g^{0}(\blam^{0})=g(\blam^{0}).
\]
 Therefore the update of our guess $\blam^{0}$ to $\blam^{1}$ is
guaranteed never to decrease the dual objective $g$. In fact, by
strict convexity of $g^{0}$, the central inequality $g^{0}(\blam^{1})\geq g^{0}(\blam^{0})$
holds with equality if and only if $\blam^{1}=\blam^{0}$, which by
assumption only holds when $\blam^{0}=\blam^{\star}$. Therefore in
fact the dual objective is guaranteed to \emph{strictly }decrease
unless $\blam^{0}=\blam^{\star}$.

The natural minorizer $g^{0}$ achieving these conditions is defined
simply as 
\[
g^{0}(\blam):=\mathbf{b}\cdot\blam-\beta^{-1}\Tr\left[e^{\beta\mathrm{diag}(\blam-\blam^{0})}e^{-\beta(C-\mathrm{diag}(\blam^{0}))}\right].
\]
 The equality $g^{0}(\blam^{0})=g(\blam)$ is immediate, and the fact
that $g^{0}\leq g$ follows directly from the \emph{Golden-Thompson
inequality} \cite{Golden1965,Thompson1965}. Strict concavity is deduced
by inspection. Provided that (\ref{eq:minor}) holds, then from our
formula (\ref{eq:blam1}) for $\blam^{1}$, it follows that $\underset{\blam}{\text{argmax}}\ g^{0}=\blam^{0}$
if and only if $\Delta\blam=0$, which following (\ref{eq:deltablam})
holds if and only if $\mathrm{diag}\left[X(\blam^{0})\right]=\mathbf{b}$,
i.e., if and only if $\blam^{0}=\blam^{\star}$.

Hence it remains only to verify (\ref{eq:minor}), i.e., that exact
maximization of the minorizer recovers our update $\blam^{0}\ra\blam^{1}$.
To this end, simplify the new trace as 
\begin{align*}
\Tr\left[e^{\beta\mathrm{diag}(\blam-\blam^{0})}e^{-\beta(C-\mathrm{diag}(\blam^{0}))}\right] & =e^{\beta(\blam-\blam^{0})}\cdot\mathrm{diag}\left[e^{-\beta(C-\mathrm{diag}(\blam^{0}))}\right]\\
 & =e^{\beta(\blam-\blam^{0})}\cdot\mathrm{diag}\left[X(\blam^{0})\right],
\end{align*}
 where $e^{\mathbf{c}}$ indicates an entrywise exponential for vectors
$\mathbf{c}$.

Minimizing $g^{0}$ can be achieved by solving the first-order optimality
condition 
\[
0=\nabla g^{0}(\blam)=\mathbf{b}-e^{\beta(\blam-\blam^{0})}\odot\mathrm{diag}\left[X(\blam^{0})\right],
\]
 where $\odot$ indicates the entrywise product of vectors. The optimality
condition is solved by isolating $\blam$ to obtain
\[
\blam=\blam^{0}+\beta^{-1}\left(\log\mathbf{b}-\log\mathrm{diag}\left[X(\blam^{0})\right]\right),
\]
 which is precisely $\blam^{1}$ as defined in (\ref{eq:blam1}).
Thus (\ref{eq:minor}) is verified.

Hence our optimization algorithm for (\ref{eq:sdpdiag}) consists
of the looping the following steps, given an initial guess for $\blam$
and a batch size $N$: 
\begin{enumerate}
\item Draw $Z\in\R^{n\times N}$ with i.i.d. $\mathcal{N}(0,1)$ entries.
\item Set $V\leftarrow Y(\mu)Z$.
\item Estimate $\mathbf{a}\approx\mathrm{diag}\left[X(\blam)\right]$ as
$\mathbf{a}\leftarrow\frac{1}{N}(V\odot V)\mathbf{1}_{N}$, following
(\ref{eq:diagest}). (Note that by construction always $\mathbf{a}>0$.)
\item Update $\blam\leftarrow\blam+\beta^{-1}\left(\log\mathbf{b}-\log\mathbf{a}\right)$.
\end{enumerate}

\subsection{Extremal eigenvalue problem: Newton's method}

The extremal eigenvalue problem (\ref{eq:sdpeig}) involves only a
scalar dual variable $\mu$, and the derivatives $g'(\mu)$ and $g''(\mu)$
are determined by (\ref{eq:eiggrad}) and (\ref{eq:eighess}). Therefore
we simply propose optimization by Newton's method.

There is considerably more flexibility available to us in the estimation
of the traces appearing in (\ref{eq:eiggrad}) and (\ref{eq:eighess}),
since we do not need to construct an entire vector of traces of size
scaling with $n$.

For thematic consistency, we will rely still on matvecs by $Y(\mu)$
defined as in (\ref{eq:Ybin}). First rewrite 
\[
g'(\mu)=k-\Tr[Y(\mu)Y(\mu)],\quad g''(\mu)=-\beta\left(\Tr[Y(\mu)Y(\mu)]-\Tr\left[Y(\mu)Y(\mu)Y(\mu)Y(\mu)\right]\right).
\]
 We can estimate both expressions using Hutchinson-type trace estimators
using only matvecs by $Y(\mu)$. Moreover our estimator for $g''(\mu)$
will preserve the negativity $g''(\mu)<0$ which is guaranteed by
concavity.

The full optimization algorithm for (\ref{eq:sdpeig}) consists of
looping the following steps, given an initial guess for $\mu$ and
a batch size $N$: 
\begin{enumerate}
\item Draw $Z\in\R^{n\times N}$ with i.i.d. $\mathcal{N}(0,1)$ entries.
\item Set $V\leftarrow Y(\mu)Z$.
\item Set $W\leftarrow Y(\mu)V$.
\item Estimate $a_{1}\approx g'(\mu)$ as $a_{1}\leftarrow\frac{1}{N}\mathbf{1}_{n}^{\top}\left(V\odot V\right)\mathbf{1}_{N}$.
\item Estimate $a_{2}\approx g''(\mu)$ as $a_{2}\leftarrow-\beta\left(a_{1}-\frac{1}{N}\mathbf{1}_{n}^{\top}(W\odot W)\mathbf{1}_{N}\right)$.
\item Update $\mu\leftarrow\mu-(k-a_{1})/a_{2}$.
\end{enumerate}

\section{Applications \label{sec:Applications}}

We consider applications on two problem types.

The first is the SDP of diagonal type (\ref{eq:sdpdiag}) arising
from the Goemans-Williamson relaxation of the Max-Cut problem. For
this problem, we show the dual solution of the entropically regularized
problem can be used to recover upper and lower bounds for the original
combinatorial problem of interest. For fixed regularization parameter
$\beta$ and assuming graphs of bounded degree, the empirical cost
scaling of solving the regularized dual SDP is only $O(n)$. Moreover,
when $\beta$ is fixed, the algorithm achieves a fixed (nontrivial)
approximation ratio as $n$ becomes large. To our knowledge, there
is no alternative $O(n)$ algorithm that can achieve such a fixed
nontrivial approximation ratio. Although our demonstration of this
scaling is only empirical, it suggests a possibility for further theoretical
analysis.

The second is the spectral embedding of graphs, which involves an
extremal eigenvalue problem for a graph Laplacian that can be rephrased
as an SDP of type (\ref{eq:sdpeig}). Usually spectral embedding is
performed by using the lowest $k$ eigenvectors to embed the graph
into $\R^{k}$. Obtaining a dual solution of the entropically regularized
problem does not give direct access to these extremal eigenvectors,
but it does give us access to matvecs by a smoothed proxy for the
spectral projector onto their span, which is sufficient to approximate
the $k$-dimensional embedding within a somewhat enlarged space of
dimension $\tilde{k}=O(k\log n)$. Here the criterion for approximating
the embedding is that the pairwise distances of the original spectral
embedding are preserved approximately. Downstream tasks such as clustering
can then be performed in the embedding space.

For graphs of bounded degree and a fixed regularization parameter
$\beta$, the total empirical cost of the dual optimization is only
$O(n)$. The subsequent recovery of the approximate spectral embedding
introduces an additional cost of $O(\tilde{k}n)$, though we comment
that the factor of $\tilde{k}$ is fully parallelizable. This scaling
should be contrasted with the $O(k^{2}n)$ scaling of direct computation
of the lowest $k$ eigenvectors by methods such as LOBPCG \cite{Knyazev2001}.
Note that since spectral embedding is a heuristic approach anyway,
it is not necessarily the case that we must approximate the spectral
embedding very accurately to reproduce its qualitative features.

\subsection{Max-Cut problem}

The Max-Cut problem \cite{GoemansWilliamson} is a combinatorial optimization
problem which can be phrased as 
\begin{equation}
\underset{x\in\{\pm1\}^{n}}{\text{minimize}}\quad x^{\top}Cx,\label{eq:maxcut}
\end{equation}
 where $C$ is usually the adjacency matrix $A$ of a graph. In fact
we shall take $C=A/n$ so that the optimal value in our experiments
remains bounded in $n$. To define related problems of Max-Cut type,
the matrix $C$ can more generally be any matrix of the same sparsity
pattern, as in the specification of spin-glass models on a graph.

The Goemans-Williamson (GW) relaxation \cite{GoemansWilliamson} of
the Max-Cut problem considers the matrix $X:=xx^{\top}$, which must
satisfy the diagonal constraint $\mathrm{diag}(X)=\mathbf{1}_{n}$
and the PSD constraint $X\succeq0$. By optimizing over $X$ and enforcing
only these conditions, we obtain a relaxation of the original problem:
\begin{align}
\underset{X\in\mathbb{R}^{n\times n}}{\text{minimize}}\quad & \Tr[CX]\label{eq:gw}\\
\text{subject to}\quad & X\succeq0,\nonumber \\
 & \mathrm{diag}(X)=\mathbf{1},\nonumber 
\end{align}
 which is an SDP of type (\ref{eq:sdpdiag}).

The dual problem to (\ref{eq:gw}) is: 
\begin{align}
\underset{\blam\in\mathbb{R}^{n}}{\text{maximize}}\quad & \mathbf{1}_{n}\cdot\blam\nonumber \\
\text{subject to}\quad & C-\mathrm{diag}(\blam)\succeq0,\label{eq:gwdual}
\end{align}
 while our \emph{regularized} dual program is 
\begin{align}
\underset{\blam\in\mathbb{R}^{n}}{\text{maximize}}\quad & \mathbf{1}_{n}\cdot\blam-\beta^{-1}\Tr\left[e^{-\beta(C-\mathrm{diag}(\blam))}\right].\label{eq:gwdualreg}
\end{align}

\subsubsection{Lower bound}

The optimal value of the relaxation (\ref{eq:gw}) furnishes a lower
bound for the Max-Cut problem (\ref{eq:maxcut}). 

Unfortunately, a dual solution $\blam^{\star}$ to (\ref{eq:gwdualreg})
does not yield a feasible point for (\ref{eq:gwdual}), which would
furnish a lower bound on the optimal value of (\ref{eq:gw}) and hence
of the Max-Cut problem (\ref{eq:maxcut}). However, for $\beta$ large,
$C-\mathrm{diag}(\blam)$ should be nearly PSD and only a small shift
should be necessary to obtain a feasible point for (\ref{eq:gwdual}).

To determine the size of the shift we need, we can compute using a
matrix-free method (which we choose to be LOBPCG \cite{Knyazev2001})
the minimal eigenvalue $\mu$ of $C-\mathrm{diag}(\blam^{\star})$.
Then $C-\mathrm{diag}(\blam^{\star}+\mu\mathbf{1}_{n})\succeq0$ is
guaranteed, i.e., $\blam^{\star}+\mu\mathbf{1}_{n}$ is feasible for
the unregularized dual problem (\ref{eq:gwdual}), and by plugging
into the dual objective we see that 
\begin{equation}
\mathbf{1}_{n}\cdot\blam^{\star}+\mu n\label{eq:lowerbound}
\end{equation}
 furnishes a lower bound on the optimal value of the Max-Cut problem
(\ref{eq:maxcut}).

\subsubsection{Upper bound}

To obtain an upper bound, a randomized rounding procedure \cite{GoemansWilliamson}
is available given the solution $X$ of the SDP relaxation (\ref{eq:gw}): 
\begin{enumerate}
\item Factorize $X=YY^{\top}$.
\item Draw $z\sim\mathcal{N}(0,I_{n})$.
\item Set $y=Yz$.
\item Set $x=\mathrm{sign}(y)\in\{\pm1\}^{n}$, where `sign' indicates the
entrywise sign.
\item Compute $x^{\top}Cx$ to obtain an upper bound on the optimal value
of the Max-Cut problem (\ref{eq:maxcut}).
\end{enumerate}
Usually the Cholesky factorization of $X$ is used in step (1), but
due to the unitary invariance of the normal distribution any choice
of factorization $X=YY^{\top}$ yields equivalent results.

We can run the same rounding procedure using the solution $X=X(\blam^{\star})$
of the \emph{regularized} problem (\ref{eq:gwdualreg}), for which
the square root factorization $Y=X^{1/2}$ is the most natural computational
choice. Indeed, note that $Y$ can be formed from our regularized
dual solution $\blam^{\star}$ as 
\[
Y=Y(\blam^{\star})=e^{-\beta(C-\mathrm{diag}(\blam^{\star}))/2}.
\]
 Note that we do not need to form the entire matrix $Y$ but instead
only require matvecs by $Y$, which are available to us by the same
procedure used within the optimization itself. In practice, we can
also repeat the rounding procedure several times and pick the lowest
upper bound recovered by the procedure, but in our experiments we
will compute the empirical expected value of the rounding procedure
over several attempts.

Using our upper and lower bounds we can construct an approximation
ratio as the ratio of the lower to upper bound.

\subsubsection{Experiments}

We consider the Max-Cut problem on Erd\H{o}s-Renyi random graphs
\cite{ErdosRenyi} of size $n$ with probability $p=3/n$ of including
an edge. Thus the expected degree of each vertex is 3. For such a
model, the rank of the optimal solution $X$ of (\ref{eq:gw}) grows
with the problem size, meaning that low-rank approaches to SDP do
not achieve $O(n)$ scaling. We will apply the noncommutative matrix
scaling algorithm of Section \ref{sec:scaling} with batch size $N=8$.

\begin{figure}[H]
\begin{centering}
\includegraphics[scale=0.6]{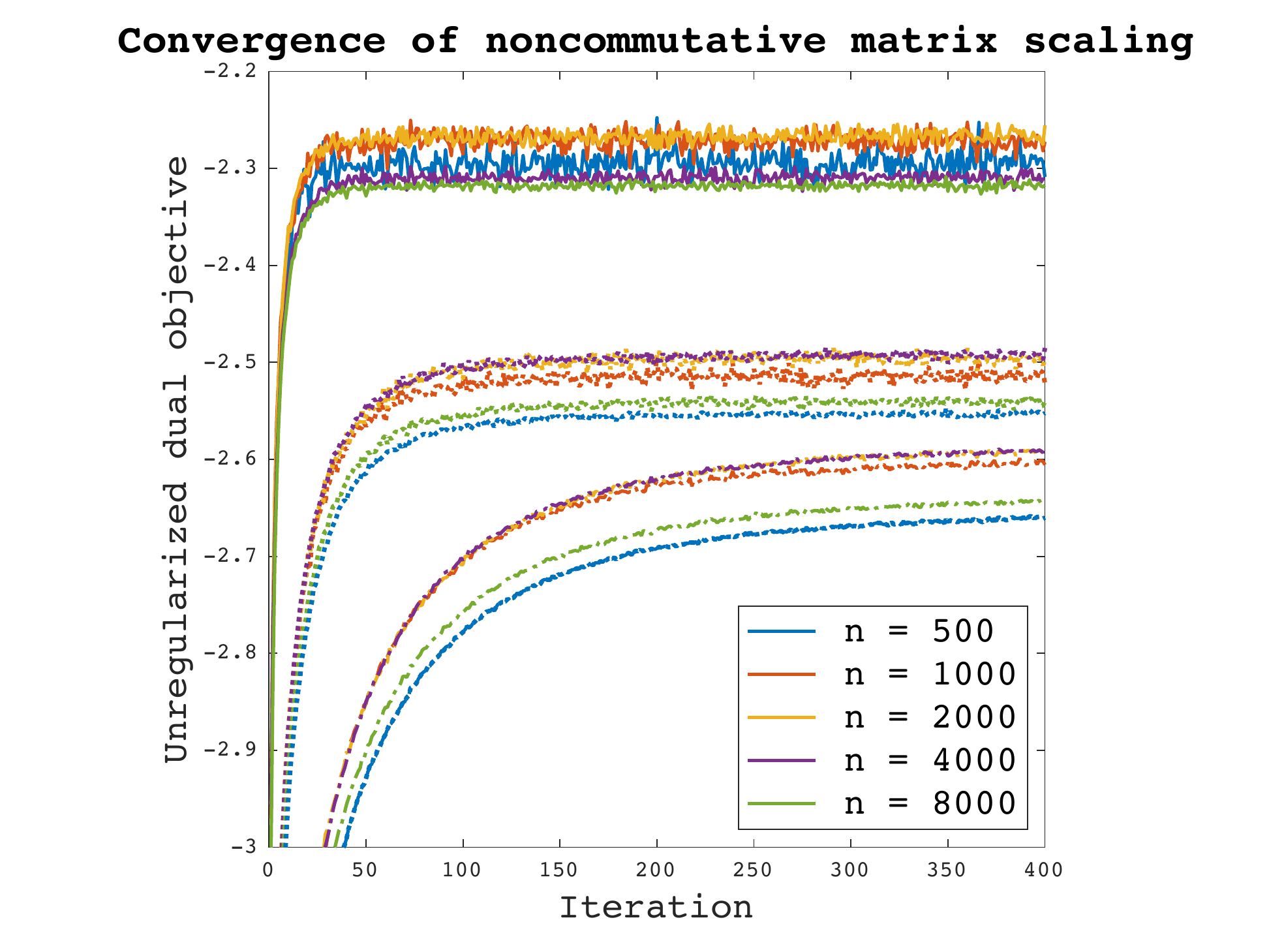}
\par\end{centering}
\caption{Convergence profile of noncommutative matrix scaling algorithm. The
solid, dotted, and dot-dashed lines indicate $\beta=10$, $\beta=32$,
and $\beta=100$, respectively, for a variety of problem sizes $n$.
The convergence rate is observed empirically to be independent of
$n$.\label{fig:maxcut_convergence}}

\end{figure}

We consider experiments in which the problem size $n$ and the regularization
parameter $\beta$ are varied. In Figure \ref{fig:maxcut_convergence},
we plot the convergence profile of the \emph{unregularized} dual objective
$\mathbf{1}_{n}\cdot\blam$. (We omit the dual regularization term
of (\ref{eq:gwdualreg}) from our convergence plots since it must
be estimated stochastically and therefore introduces spurious noise
in the plots.) We see that for fixed $\beta$, the convergence rate
is empirically independent of the problem size $n$. We also see that
convergence is slower when $\beta$ is large, though larger $\beta$
allows us to drive the unregularized dual objective down. (Note that
in addition the cost of matvecs by $Y(\blam)$ scales linearly with
$\beta$.)

When $\beta$ is small, the regularization term might change the original
unregularized SDP (\ref{eq:gw}) significantly. Therefore we aim to
quantify the impact of $\beta$ on the approximation ratio of the
problem. It is known \cite{GoemansWilliamson} that the GW relaxation
achieves an approximation ratio of 
\[
\alpha=\frac{2}{\pi}\min_{\theta\in[0,\pi]}\frac{\theta}{1-\theta}\approx0.878,
\]
 which is in fact the best possible approximation ratio if the unique
games conjecture is true \cite{UniqueGames}. Note that a trivial
approximation ratio of $1/2$ is always available (in expectation)
by drawing each of the entries of $x$ independently randomly from
the uniform distribution over $\{\pm1\}$.

\begin{figure}[H]
\begin{centering}
\includegraphics[scale=0.6]{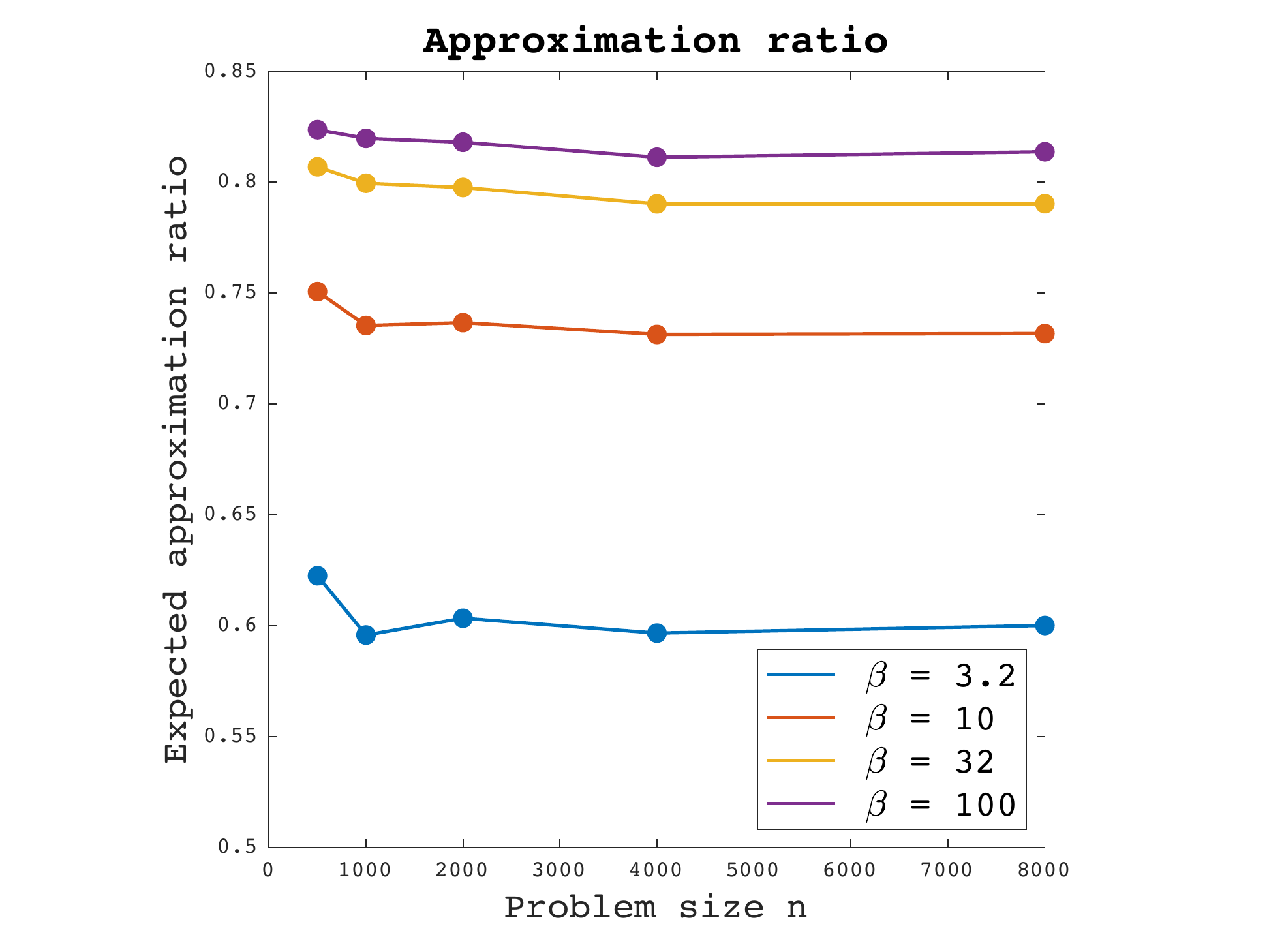}
\par\end{centering}
\caption{Approximation ratio of our algorithm as a function of problem size
$n$ for various values of $\beta$. We observe empirically that a
constant nontrivial approximation ratio is maintained as $n$ becomes
large. The ratio is improved as $\beta$ increases. \label{fig:maxcut_ratio}}
\end{figure}

In Figure \ref{fig:maxcut_ratio}, we consider the effect of fixing
$\beta$ and increasing the problem size in order to determine whether
a constant approximation ratio is achieved as $n\ra\infty$. We compute
the approximation ratio as the ratio of the lower bound furnished
as specified above to the\emph{ expected} lower bound (computed as
an empirical average over 1000 samples) furnished as specified by
steps 1-5 above. The optimization algorithm is run with batch size
$N=8$ for 400 iterations in each case, except the case $\beta=100$
in which the optimization is run over 1000 iterations to ensure convergence.
In the figure, we see that a constant nontrivial approximation ratio
is maintained as $n$ becomes large. Moreover, the ratio is improved
by enlarging $\beta$.

\subsection{Spectral embedding}

The task of spectral embedding \cite{PothenEtAl1990,SpectralClustering}
is to embed a graph into a Euclidean space in which the Euclidean
distance is meaningful. Oftentimes, the graph is constructed from
a point cloud in a preprocessing step via the selection of nearest
neighbors or a procedure making use of a similarity kernel. After
embedding, downstream tasks such as clustering can be significantly
easier.

There are several approaches to spectral clustering given a graph,
but a standard approach \cite{SpectralClustering} is specified as
follows. Given the adjacency matrix $A$ of the graph, let $D=\mathrm{diag}(A\mathbf{1}_{n})$
denote the diagonal matrix with diagonal entries given by the degrees
of the corresponding graph vertices. Then a symmetric normalized graph
Laplacian is constructed as 
\[
L=I-D^{-1/2}AD^{-1/2}.
\]

If we desire a spectral embedding into $\R^{k}$, orthonormal eigenvectors
corresponding to the $k$ lowest eigenvectors of $L$ are collected
into the columns of $\Phi\in\R^{n\times k}$. The map $i\mapsto\Phi_{i,:}\in\R^{k}$
to the $i$-th row of matrix $\Phi$ then defines the spectral embedding
map.

We rephrase the problem of determining $\Phi$ into the problem of
determining $X=\Phi\Phi^{\top}$ via the SDP (\ref{eq:sdpeig}) in
which we take $C=L$, and we introduce the binary von Neumann entropy
with finite temperature $\beta$ as a regularization term.

\subsubsection{Recovering the spectral embedding \label{sec:recovery}}

Suppose we are given the regularized dual solution $\mu^{\star}$
and for simplicity denote $X=X(\mu^{\star})$ and $Y=Y(\mu^{\star})$.
We now describe in principle how to recover an approximate spectral
embedding from $\mu^{\star}$.

First create a random matrix $Z\in\R^{n\times\tilde{k}}$ consisting
of entries drawn independently from the standard normal distribution.
We denote the $i$-th column of $Z$ as $z^{(i)}$ and form 
\[
\Psi:=\frac{1}{\sqrt{\tilde{k}}}YZ\in\R^{n\times\tilde{k}}
\]
 with $i$-th column $\psi^{(i)}=\frac{1}{\sqrt{\tilde{k}}}Yz^{(i)}$.
Note that forming $\Psi$ only requires matvecs by $Y=Y(\mu^{\star})=F_{\beta}^{1/2}(C-\mu^{\star}I)$,
as in the optimization procedure itself.

We claim that the rows of $\Psi$ furnish an approximate spectral
embedding in the sense that the Gram matrix of inner products $G:=\Psi\Psi^{\top}$
well approximates $X(\mu^{\star})=F_{\beta}(C-\mu I)$, which is itself
a smoothed approximation to the spectral projector $\Phi\Phi^{\top}$
onto the lowest $k$ eigenvectors of $C$. Since pairwise inner products
determine all the pairwise Euclidean distances, this implies that
the embeddings $\Psi_{i,:}\in\R^{\tilde{k}}$ furnish an effective
approximate solution to the original problem.

We want to justify that we can take $\tilde{k}=O(k\log n)$ in order
to achieve a good approximation $\Psi\Psi^{\top}\approx X(\mu^{\star}).$
Then, assuming that the optimization furnishing $\mu^{\star}$ can
be converged in $O(1)$ iterations, it would follow that approximate
spectral embedding can be achieved in $O(kn\log n)$ time, which for
large values of $k$ should outperform the usual $O(k^{2}n)$ scaling
of determining $\Phi$ exactly.
\begin{prop}
Let $\delta\in(0,1/2]$, and suppose that $z^{(i)},i=1,\ldots,N$
are i.i.d. with entries that are themselves i.i.d. sub-Gaussian random
variables with a constant sub-Gaussian parameter, mean zero, and unit
variance. There exist constants $c,C>0$ such that if $\tilde{k}>c\log\left(\frac{n}{\delta}\right)$,
then with probability at least $1-\delta$, the inequality 
\[
\vert G_{ij}-X_{ij}\vert\leq C\sqrt{\log\left(\frac{n}{\delta}\right)}\ \sqrt{\frac{X_{ii}X_{jj}}{\tilde{k}}}
\]
 holds for all $i,j=1,\ldots,n$. It follows that under the same conditions,
\[
\frac{\Vert G-X\Vert_{\mathrm{F}}}{\Vert X\Vert_{\mathrm{F}}}\leq C\sqrt{\log\left(\frac{n}{\delta}\right)}\ \sqrt{\frac{k}{\tilde{k}}}\ \frac{\sqrt{k}}{\Vert X\Vert_{\mathrm{F}}}.
\]
\end{prop}

\begin{rem}
Note that if $X=\Phi\Phi^{\top}$, then $\Vert X\Vert_{\mathrm{F}}=\sqrt{k}$.
Even for finite $\beta$ where $X\approx\Phi\Phi^{\top}$ holds only
approximately, for most reasonable models of the eigenvalues of the
graph Laplacian $L$, we can expect that $\frac{\sqrt{k}}{\Vert X\Vert_{\mathrm{F}}}=O(1)$
for fixed $\beta$. Under this condition, to obtain an approximate
spectral embedding of some fixed relative accuracy (in the sense of
relative Frobenius norm error of the Gram matrix), we need only take
$\tilde{k}=O(k\log n)$.
\end{rem}

\begin{proof}
Note that 
\[
G=\Psi\Psi^{\top}=\sum_{l=1}^{\tilde{k}}\psi^{(l)}\psi^{(l)\top}=\frac{1}{\tilde{k}}\sum_{l=1}^{\tilde{k}}Yz^{(l)}z^{(l)\top}Y.
\]
Hence we can think of $G$ as an empirical expectation with $\E[G]=X$.
In fact we can even think of it as a Hutchinson trace estimator for
each entry $G_{ij}$ by rearranging 
\[
G_{ij}=e_{i}^{\top}Ge_{j}=\frac{1}{\tilde{k}}\sum_{l=1}^{\tilde{k}}z^{(l)\top}\left[Ye_{j}e_{i}^{\top}Y\right]z^{(l)}.
\]
 Hence $G_{ij}$ is a Hutchinson estimator for the trace $\Tr[Ye_{j}e_{i}^{\top}Y]$.
Since 
\[
\Vert Ye_{j}e_{i}^{\top}Y\Vert_{\mathrm{F}}=\Tr\left[Ye_{j}e_{i}^{\top}YYe_{i}e_{j}^{\top}Y\right]=X_{ii}X_{jj},
\]
 we can again directly apply Lemma 2.1 of \cite{Hutch++} to each
Hutchinson estimator, $i,j=1,\ldots,n$, which, together with the
union bound over the $n^{2}$ separate traces, implies the first statement
of the proposition.

The second statement follows by recalling that $X=X(\mu^{\star})$
satisfies the constraint $\Tr[X]=k$, since it must be primal feasible.
\end{proof}

\subsubsection{Experiments}

We consider a random graph model with $n$ vertices in which $n$
is a multiple of $m$ and the subgraphs corresponding to the vertex
subsets $\{1,\ldots,m\}$, $\{m+1,\ldots,2m\}$, ..., and $\{n-m+1,\ldots,n\}$
are fully connected. In addition, every edge connecting two vertices
appearing in distinct subsets is included with probability $1/n$.
Thus the graph should be viewed as having $k=n/m$ clusters, which
should be captured well by spectral embedding into $\R^{k}$.

As the approximation $X(\mu^{\star})=\Phi\Phi^{\top}$ can be well-understood
in terms of the spectrum of $L$ and as the theory of recovering the
spectral embedding is well-understood following Section \ref{sec:recovery},
we focus in our experiments on the optimization, which is the point
of more general interest for SDP with upper and lower semidefinite
constraints.

Since the eigenvalues of the symmetric normalized graph Laplacian
are bounded between $0$ and $2$ \cite{SpectralGraphTheory}, following
Section \ref{sec:matvecs}, we can take $R=2$ for our polynomial
approximation of $F_{\beta}^{1/2}$. We insist on a sup norm error
of at most $10^{-5}$ for our polynomial approximation over the interval
$[-2,2]$.

\begin{figure}[H]
\begin{centering}
\includegraphics[scale=0.6]{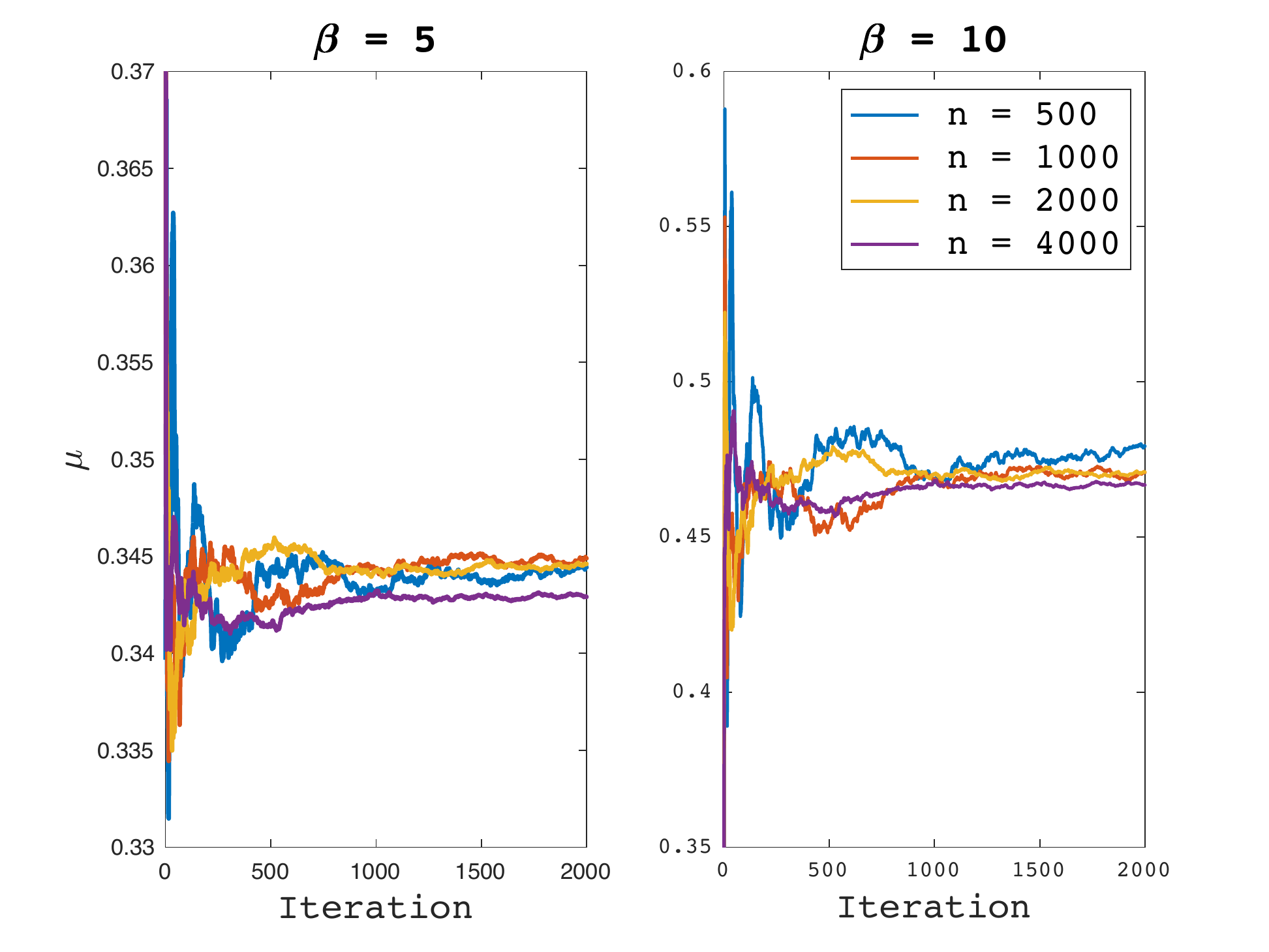}
\par\end{centering}
\caption{Smoothed optimization trajectory of $\mu$ for various problem sizes
$n$. Left: $\beta=5$. Right: $\beta=10$. The convergence rate is
observed to be empirically independent of $n$.\label{fig:chebopt}}
\end{figure}

In our experiments we take $m=10$, so $k=n/10$. Moreover, instead
of computing each Newton step with high accuracy, we instead opt to
take a small batch size $N=8$. Due to the large noise in the update,
we can view our iterative optimization method as a stochastic process.
At the $t$-th iteration, we compute a running average of the values
of $\mu$ over the preceding $\left\lfloor t/2\right\rfloor $ iterations
to obtain a convergent trajectory. In all experiments we run the procedure
for $2000$ total iterations, and we plot the smoothed trajectories
in Figure \ref{fig:chebopt}. At the end of the optimization, we verified
using a Hutchinson estimator with $1000$ random vectors that the
primal trace constraint $\Tr[X(\mu^{\star})]=k$ is satisfied to within
relative error of 0.01 in all of our experiments. From Figure \ref{fig:chebopt}
we observe empirically that the convergence rate is independent of
$n$.

\bibliographystyle{plain}
\bibliography{bigbib}

\end{document}